\documentclass[a4paper, 12pt]{article}

\usepackage{fullpage}
\usepackage{authblk}
\usepackage{amsmath}
\usepackage{amssymb}
\usepackage{amsthm}
\usepackage{xfrac}
\usepackage{url}
\usepackage{csquotes}
\usepackage{enumitem}

\MakeOuterQuote{"}
\setlist{nolistsep}
\setlist[itemize]{leftmargin=*}
\setlist[enumerate]{leftmargin=*}

\title{Linear Algebraic Number Theory, \\ Part I: Foundations}
\author{Joram Soch}
\affil{BCCN Berlin, Germany}
\affil{joram.soch@bccn-berlin.de}
\date{}

\begin{document}

\setcounter{page}{0}
\maketitle

\begin{abstract}
\noindent
We introduce a new framework called linear algebraic number theory (LANT) that reformulates the number-theoretic problem as a regression model and solves it using matrix algebra. This framework restricts all computations to log space, therefore replaces multiplication with addition and allows to capture variation in the natural numbers from variation in the prime numbers. This automatically puts prime numbers to their designated place of atomic particles of natural numbers and enables fruitful new formulations of number-theoretic functions. We outline the theory, derive some basic results, make connections to standard number theory and give an outlook regarding the Riemann hypothesis, number theory's long-standing enigma.
\end{abstract}

\vspace{1em}
\tableofcontents

\pagebreak
\section{Introduction}

Let $n$ be a positive natural number. The \textit{fundamental theorem of arithmetic} states that there is a unique factorization by which $n$ can be written as a product of prime powers:

\begin{equation} \label{eq:FTA}
n = \prod_{i=1}^k p_i^{n_i} \; .
\end{equation}

In this product, $k$ is the number of primes that divide $n$, $p_1 < \ldots < p_k$ are prime numbers and $n_1, \ldots, n_k$ are positive integers. This is called the \textit{canonical representation} or \textit{standard form} of $n$. For example,

\begin{equation}
360 = 2^3 \times 3^2 \times 5^1 \; .
\end{equation}

We can take the natural logarithm of equation (\ref{eq:FTA}) and obtain

\begin{equation} \label{eq:FTA-ln}
\ln n = \sum_{i=1}^k n_i \ln p_i \; .
\end{equation}

Applied to the example, this gives

\begin{equation}
\ln 360 = 3 \ln 2 + 2 \ln 3 + 1 \ln 5 \; .
\end{equation}

We can write equation (\ref{eq:FTA-ln}) as a vector product and obtain

\begin{equation} \label{eq:FTA-ln-VP}
\ln n = \begin{bmatrix} n_1 & \cdots & n_k \end{bmatrix} \begin{bmatrix} \ln p_1 \\ \vdots \\ \ln p_k \end{bmatrix} \; .
\end{equation}

Applied to the example, this gives

\begin{equation}
\ln 360 = \begin{bmatrix} 3 & 2 & 1 \end{bmatrix} \begin{bmatrix} \ln 2 \\ \ln 3 \\ \ln 5 \end{bmatrix} \; .
\end{equation}

Note that, as $p^0 = 1$ for any $p \in \mathbb{R}, p \neq 0$, one might insert prime powers with exponent zero in equation (\ref{eq:FTA}) or log primes with factor zero in equation (\ref{eq:FTA-ln}) without changing the value of $n$. This means that the second vector in equation (\ref{eq:FTA-ln-VP}) will be the same for all $n$, namely a column vector of all log primes, and just the first vector has to be adapted in order to achieve the correct combination of primes.

The basic idea of this paper is to make use of this insight and rewrite the prime factorization of natural numbers (or, log-prime summation of log integers) as a linear equation system which has the logarithmized natural numbers on its left-hand side and a matrix product of a factorization matrix and the log primes on its right-hand side. This is referred to as \textit{linear algebraic number theory} (LANT).

\pagebreak
\section{Definitions}

An example for such a linear equation system is

\begin{equation}
\begin{split}
\ln 1 &= 0 \ln 2 + 0 \ln 3 + 0 \ln 5 \\
\ln 2 &= 1 \ln 2 + 0 \ln 3 + 0 \ln 5 \\
\ln 3 &= 0 \ln 2 + 1 \ln 3 + 0 \ln 5 \\
\ln 4 &= 2 \ln 2 + 0 \ln 3 + 0 \ln 5 \\
\ln 5 &= 0 \ln 2 + 0 \ln 3 + 1 \ln 5 \\
\ln 6 &= 1 \ln 2 + 1 \ln 3 + 0 \ln 5
\end{split}
\end{equation}

which, in matrix algebra notation, can be written as

\begin{equation} \label{eq:FTA-ln-LES}
\begin{bmatrix} \ln 1 \\ \ln 2 \\ \ln 3 \\ \ln 4 \\ \ln 5 \\ \ln 6 \end{bmatrix} =
\begin{bmatrix} 0 & 0 & 0 \\ 1 & 0 & 0 \\ 0 & 1 & 0 \\ 2 & 0 & 0 \\ 0 & 0 & 1 \\ 1 & 1 & 0 \end{bmatrix}
\begin{bmatrix} \ln 2 \\ \ln 3 \\ \ln 5 \end{bmatrix} \; .
\end{equation}

In order to formulate this for the general case, we will introduce some definitions.

\begin{quote}
\textbf{Definition 1:} (\textit{element-wise logarithm}) Whenever the natural logarithm is applied to a vector $v \in \mathbb{R}^n$, it is calculated element-wise:
\begin{equation} \label{eq:ln-v}
\ln v = \ln \begin{bmatrix} v_1 \\ \vdots \\ v_n \end{bmatrix} = \begin{bmatrix} \ln v_1 \\ \vdots \\ \ln v_n \end{bmatrix} \; .
\end{equation}
\end{quote}

\begin{quote}
\textbf{Definition 2:} (\textit{natural number vector}) Let $n$ be a positive natural number. Then, the $n \times 1$ vector $z_n$ is defined as
\begin{equation} \label{eq:zn}
z_n = \begin{bmatrix} 1 \\ 2 \\ \vdots \\n \end{bmatrix} \; .
\end{equation}
\end{quote}

\begin{quote}
\textbf{Definition 3:} (\textit{prime number vector}) Let $n$ be a positive natural number. Then, the $\pi(n) \times 1$ vector $p_n$ is defined as
\begin{equation} \label{eq:pn}
p_n = \begin{bmatrix} 2 \\ 3 \\ 5 \\ \vdots \\ p \end{bmatrix}
\end{equation}
where $p$ is the largest $x \in \mathbb{P}$ for which $x \leq n$, $\mathbb{P}$ is the set of prime numbers and $\pi(n)$ is the number of primes less than or equal to $n$.
\end{quote}

Obviously, as implied by the fundamental theorem of arithmetic (\ref{eq:FTA}) and instantiated by the above example (\ref{eq:FTA-ln-LES}), $z_n$ and $p_n$ can be related to each other in log space by a matrix of coefficients. We will call this the \textit{factorization matrix}.

\begin{quote}
\textbf{Definition 4:} (\textit{prime factorization matrix}) Let $n$ be a positive natural number. Then, the prime factorization matrix is the $n \times \pi(n)$ matrix $F_n$ for which
\begin{equation} \label{eq:Fn}
\ln z_n = F_n \ln p_n \; .
\end{equation}
\end{quote}

We know that $F_n$ exists for a given $n > 1$, because every natural number greater 1 is either prime or can be factorized into primes smaller than itself, and equation (\ref{eq:Fn}) is nothing but a restatement of that fact. We also know that $F_n$ is unique as there is only one prime factorization for every $n > 1$ which can be proven with recourse to Euclid's lemma (Euclid, VII, 30). An example for $n = 20$ is given in Figure~1.

\begin{center}
\includegraphics[width=0.75\linewidth, clip=true, trim=50 160 25 80]{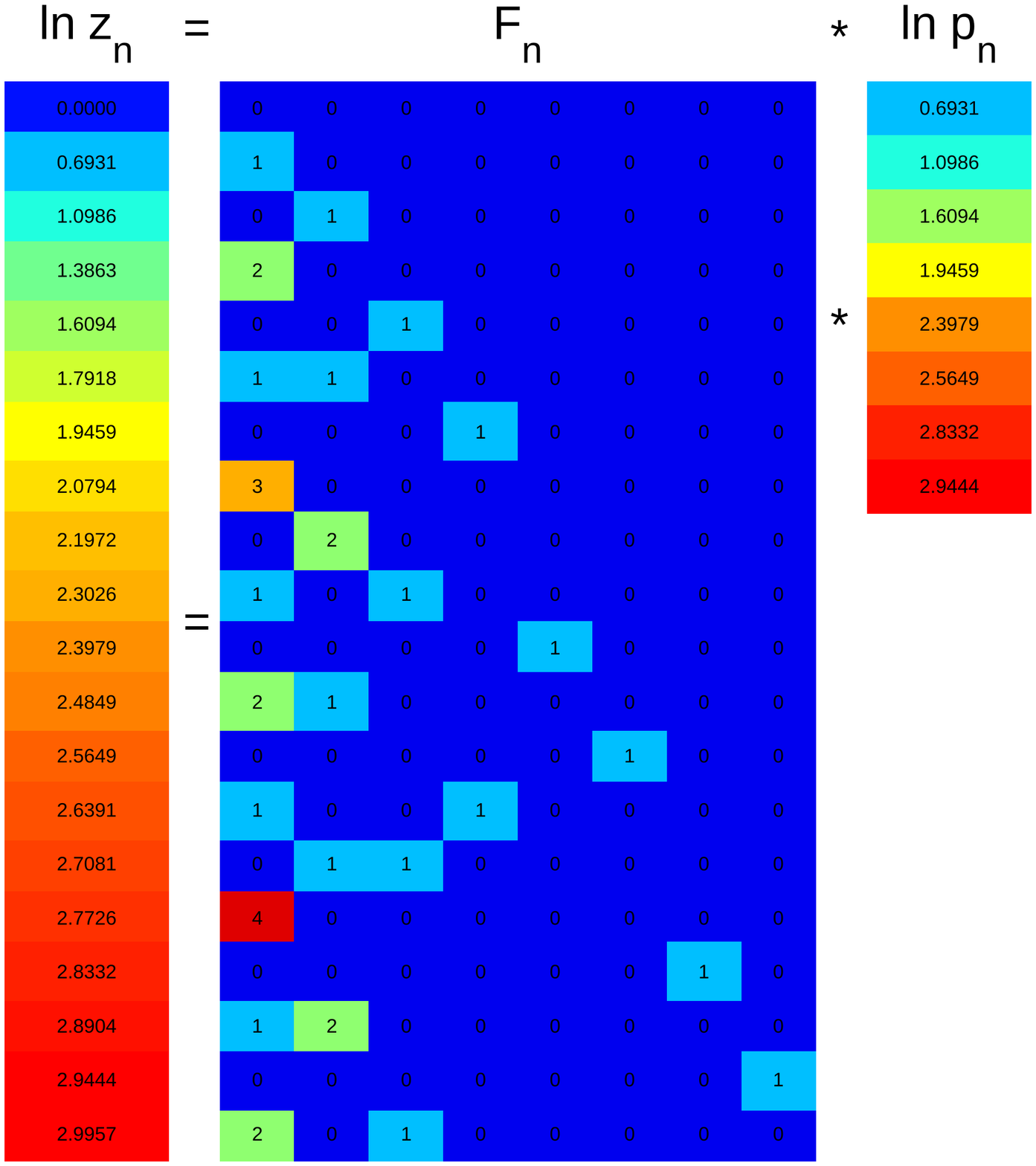}
\end{center}

\textbf{Figure 1.} Log integers $\ln z_n$, factorization matrix $F_n$ and log primes $\ln p_n$ for $n = 20$. This figure illustrates the log-space analogue of integer factorization in which log integers are represented as sums of log primes, weighted by the factorization matrix.

\pagebreak
Suppose we didn't know the primes up to a certain number $n$. A natural consequence would be that we try to derive, solve for, infer on or estimate them. To this end, we introduce the concept of \textit{candidate primes}.

\begin{quote}
\textbf{Definition 5:} (\textit{candidate prime vector}) Let $n$ be a positive natural number. A candidate prime vector $q$ is an $m \times 1 $ vector with $m \leq n$ which only contains pairwise different natural numbers smaller than or equal to $n$.
\end{quote}

For example, the following would be possible candidate primes for $n = 10$:

\begin{equation} \label{eq:q-ex}
q_1 = \begin{bmatrix} 2 \\ 5 \end{bmatrix}, \quad q_2 = \begin{bmatrix} 2 \\ 3 \\ 5 \\ 7 \end{bmatrix} = p_{10}, \quad q_3 = \begin{bmatrix} 2 \\ 3 \\ 4 \\ 5 \\ 7 \\ 8 \end{bmatrix}, \quad q_4 = \begin{bmatrix} 1 \\ 2 \\ 3 \\ 4 \\ 5 \\ 6 \\ 7 \\ 8 \\ 9 \\ 10 \end{bmatrix} = z_{10} \; .
\end{equation}

It would now be tempting to take a certain vector $q$ and somehow calculate its associated matrix $F_n$ in order to factorize potentially large numbers. However, this is neither possible nor necessary. It is not possible as this equation system would contain more unknowns than equations. It is not necessary as every set of candidate primes already implies a factorization matrix which we call a \textit{candidate factorization}.

\begin{quote}
\textbf{Definition 6:} (\textit{candidate factorization matrix}) Let $q$ be a candidate prime vector. The candidate factorization matrix $F_n(q)$ is the $n \times m$ matrix that would be the prime factorization matrix, if $q$ were the true primes.
\end{quote}

For example, as 2 is a prime number, every second number is factorized by $2^1$, every fourth number is factorized by $2^2$, every eigth number is factorized by $2^3$ and so on. Consequently, the first column of $F_n$ has a 1 in every second row, a 2 in every fourth row, a 3 in every eigth row (see Figure~1). Similarly, if 4 were prime (which it is not), the corresponding column of $F_n$ would have a 1 in every fourth row, a 2 in every sixteenth row, a 3 in every sixty-fourth row etc.

To continue with the example from above, the candidate factorizations for the candidate primes in equation (\ref{eq:q-ex}) would be:

\begin{equation} \label{eq:Fnq-ex1}
F_{10}(q_1) = \begin{pmatrix} 0 & 0 \\ 1 & 0 \\ 0 & 0 \\ 2 & 0 \\ 0 & 1 \\ 1 & 0 \\ 0 & 0 \\ 3 & 0 \\ 0 & 0 \\ 1 & 1 \\ \end{pmatrix}, \;
F_{10}(q_2) = \begin{pmatrix} 0 & 0 & 0 & 0 \\ 1 & 0 & 0 & 0 \\ 0 & 1 & 0 & 0 \\ 2 & 0 & 0 & 0 \\ 0 & 0 & 1 & 0 \\ 1 & 1 & 0 & 0 \\ 0 & 0 & 0 & 1 \\ 3 & 0 & 0 & 0 \\ 0 & 2 & 0 & 0 \\ 1 & 0 & 1 & 0 \end{pmatrix}, \;
F_{10}(q_3) = \begin{pmatrix} 0 & 0 & 0 & 0 & 0 & 0 \\ 1 & 0 & 0 & 0 & 0 & 0 \\ 0 & 1 & 0 & 0 & 0 & 0 \\ 2 & 0 & 1 & 0 & 0 & 0 \\ 0 & 0 & 0 & 1 & 0 & 0 \\ 1 & 1 & 0 & 0 & 0 & 0 \\ 0 & 0 & 0 & 0 & 1 & 0 \\ 3 & 0 & 1 & 0 & 0 & 1 \\ 0 & 2 & 0 & 0 & 0 & 0 \\ 1 & 0 & 0 & 1 & 0 & 0 \end{pmatrix},
\end{equation}

\begin{equation} \label{eq:Fnq-ex2}
F_{10}(q_4) = \begin{pmatrix} 1 & 0 & 0 & 0 & 0 & 0 & 0 & 0 & 0 & 0 \\ 1 & 1 & 0 & 0 & 0 & 0 & 0 & 0 & 0 & 0 \\ 1 & 0 & 1 & 0 & 0 & 0 & 0 & 0 & 0 & 0 \\ 1 & 2 & 0 & 1 & 0 & 0 & 0 & 0 & 0 & 0 \\ 1 & 0 & 0 & 0 & 1 & 0 & 0 & 0 & 0 & 0 \\ 1 & 1 & 1 & 0 & 0 & 1 & 0 & 0 & 0 & 0 \\ 1 & 0 & 0 & 0 & 0 & 0 & 1 & 0 & 0 & 0 \\ 1 & 3 & 0 & 1 & 0 & 0 & 0 & 1 & 0 & 0 \\ 1 & 0 & 2 & 0 & 0 & 0 & 0 & 0 & 1 & 0 \\ 1 & 1 & 0 & 0 & 1 & 0 & 0 & 0 & 0 & 1 \end{pmatrix} \; .
\end{equation}

To express general candidate factorizations, we will introduce some more definitions.

\begin{quote}
\textbf{Definition 7:} (\textit{basic vectors}) The zero vector and the ones vector:
\begin{equation} \label{eq:0n-1n}
0_n = \left. \begin{bmatrix} 0 \\ \vdots \\ 0 \end{bmatrix} \right\rbrace \text{n zeros}, \quad
1_n = \left. \begin{bmatrix} 1 \\ \vdots \\ 1 \end{bmatrix} \right\rbrace \text{n ones} \; .
\end{equation}
\end{quote}

\begin{quote}
\textbf{Definition 8:} (\textit{elementary vectors}) The $i$-th elementary vector in $n$-dimensional vector space is an $n$-dimensional zero vector with a one in its $i$-th entry:
\begin{equation} \label{eq:ev}
e_{i|n} = \begin{bmatrix} 0 \\ \vdots \\ 1 \\ \vdots \\ 0 \end{bmatrix} \begin{matrix} \vphantom{0} \\ \vphantom{\vdots} \\ \leftarrow \text{i-th position} \\ \vphantom{\vdots} \\ \leftarrow \text{n-th position} \end{matrix} \; .
\end{equation}
\end{quote}

\begin{quote}
\textbf{Definition 9:} (\textit{periodic elementary vectors}) The $i$-th periodic elementary vector in $n$-dimensional vector space is an $n$-dimensional zero vector with a one in its $k$-th entries where $k = i, 2i, \ldots, \lfloor \sfrac{n}{i} \rfloor i$:
\begin{equation} \label{eq:pev}
e_{\bar{i}|n} = \sum_{k=1}^{\lfloor \sfrac{n}{i} \rfloor} e_{k \cdot i|n} = \begin{bmatrix} 1_{\lfloor \sfrac{n}{i} \rfloor} \otimes e_{i|i} \\ 0_{\mathrm{mod}(n,i)} \end{bmatrix} = \begin{bmatrix} 0 \\ \vdots \\ 0 \\ 1 \\ 0 \\ \vdots \\ 0 \\ 1 \\ 0 \\ \vdots \\ 0 \end{bmatrix} \begin{matrix} \vphantom{0} \\ \vphantom{\vdots} \\ \vphantom{0} \\ \leftarrow \text{i-th position} \\ \vphantom{0} \\ \vphantom{\vdots} \\ \vphantom{0} \\ \leftarrow \text{2i-th position} \\ \vphantom{0} \\ \vphantom{\vdots} \\ \leftarrow \text{n-th position} \end{matrix} \; .
\end{equation}
\end{quote}

\pagebreak
With these definitions, we are now able to express arbitrary factorization matrices:

\begin{quote}
\textbf{Definition 10:} (\textit{factorization vector}) The $i$-th factorization vector in $n$-dimensional vector space describes how often $i$ would occur as a factor in the prime factorization of the numbers $z_n$, if $i$ were prime. It is given by
\begin{equation} \label{eq:fv}
f_{i|n} = \sum_{j=1}^{\lfloor \log_{i}n \rfloor} e_{\overline{i^j}|n} \; .
\end{equation}
By definition, we set
\begin{equation} \label{eq:f1}
f_{1|n} = 1_n \; .
\end{equation}
\end{quote}

\begin{quote}
\textbf{Theorem 1:} Let $q$ be a candidate prime vector. Then, the corresponding candidate factorization matrix is given by
\begin{equation} \label{eq:Fnq}
F_n(q) = \left[ f_{q_1|n} \; \ldots \; f_{q_m|n} \right] \; .
\end{equation}
\end{quote}

\begin{quote}
\textbf{Proof 1:} This follows from Def.~6 and 10. The sum over periodic elementary vectors in equation (\ref{eq:fv}) ensures that the $j$-th power of each candidate prime $q_i$ repeats every $q_i^j$-th entry, because $q_i^j$ would be part of these factorizations, if $q_i$ was prime. \hspace\fill $\blacksquare$
\end{quote}

\pagebreak
We conclude this section with two definitions that will later become relevant for the inversion of quadratic candidate factorization matrices.

\begin{quote}
\textbf{Definition 11:} (\textit{elementary matrix}) An elementary matrix is an $n \times n$ matrix which performs an elementary row operation on another $n \times n$ matrix when being multiplied from the left:
\begin{equation} \label{eq:EM-Pij}
P_{ij} = \begin{bmatrix} 1 \\ & \ddots \\ & & 0 & & 1 \\ & & & \ddots \\ & & 1 & & 0 \\ & & & & & \ddots \\ & & & & & & 1 \end{bmatrix} \; ,
\end{equation}
\begin{equation} \label{eq:EM-Mi}
M_i(\lambda) = \begin{bmatrix} 1 \\ & \ddots \\ & & 1 \\ & & & \lambda \\ & & & & 1 \\ & & & & & \ddots \\ & & & & & & 1 \end{bmatrix} \; ,
\end{equation}
\begin{equation} \label{eq:EM-Gij}
G_{ij}(\lambda) = \begin{bmatrix} 1 \\ & \ddots \\ & & 1 \\ & & & \ddots \\ & & \lambda & & 1 \\ & & & & & \ddots \\ & & & & & & 1 \end{bmatrix} \; .
\end{equation}
$P_{ij}$ exchanges rows $i$ and $j$ of a matrix, $M_i(\lambda)$ multiplies the $i$-th row of a matrix with $\lambda$ and $G_{ij}(\lambda)$ multiplies row $j$ by $\lambda$ and adds it to row $i$ of a matrix when being multiplied from the left.
\end{quote}

\begin{quote}
\textbf{Definition 12:} (\textit{extended elementary matrix}) An extended elementary matrix is a matrix which performs more than one row operation when being multiplied from the left. An example would be
\begin{equation} \label{eq:EEM}
A = M_3(8) + M_7(5) - I_{10}
\end{equation}
which multiplies the third row with 8 and the seventh row with 5.
\end{quote}

\pagebreak
\section{Modelling the Natural Numbers}

The factorization matrix is at the heart of LANT and appears in its fundamental theorem:

\begin{equation} \label{eq:FT-LANT}
\ln z_n = F_n \ln p_n \; .
\end{equation}

We will now again assume that $p_n$ is unknown, so that we have to solve for it. This can be nicely connected to linear models and statistical modelling as, when searching the optimal solution for $\ln p_n$, we are seeking the best way to capture the values of the natural numbers, just like we are seeking the best way to capture the variance in measured data when applying linear models to empirical phenomena.

The univariate linear regression model is given by

\begin{equation} \label{eq:GLM}
y = X \beta + \varepsilon \; .
\end{equation}

In this equation, certain data ($n \times 1$ vector $y$) are modelled as a linear combination of independent variables ($n \times p$ matrix $X$), weighted by some coefficients ($p \times 1$ vector $\beta$), plus some residuals that cannot be explained ($n \times 1$ vector $\varepsilon$) where $y$ is called the \textit{signal}, $X$ is called the \textit{design matrix}, $\beta$ are called \textit{regression coefficients} and $\varepsilon$ is called \textit{noise}. We observe the following parallels between (\ref{eq:FT-LANT}) and (\ref{eq:GLM}):

\vspace{0.5em}
\begin{itemize}
	
\setlength{\itemsep}{0em}
\setlength{\parskip}{0.5em}
	
\item
The log integers $\ln z_n$ are the signal $y$ that we want to explain.

\item
The factorization matrix $F_n$ is the design matrix $X$ that we use to explain.

\item
The log primes $\ln p_n$ are the regression coefficients $\beta$ that we want to estimate.

\item
If we use the prime factorization matrix $F_n$ as the design matrix, there are no residuals $\varepsilon$ as the natural numbers are completely explained by the prime numbers (see Figure~1). However, if we use a candidate factorization matrix lacking some primes, we will fail to resolve the complete variation in the natural numbers, so that there will be errors $\varepsilon$ (see Figure~2).
	
\end{itemize}
\vspace{0.5em}

We can therefore write down the statistical version of equation (\ref{eq:FT-LANT}):

\begin{equation} \label{eq:LANT-GLM}
\ln z_n = F_n(q) \ln q + \varepsilon_n \; .
\end{equation}

With the concepts of candidate primes and candidate factorization, we have a simple method of constructing the design matrix for our linear regression (\ref{eq:LANT-GLM}). The next step is therefore to estimate the model, i.e. to find some parameters, given the data and the design:

\begin{equation} \label{eq:GLM-est}
\hat{\beta} = f(y,X) \; .
\end{equation}

Naturally, when performing linear regression, one wants to keep the residuals $\varepsilon$ as small as possible in order to achieve "the best possible fit" of the model to the data. A common framework for assigning parameter values following this rationale is \textit{ordinary least squares} (OLS).

\pagebreak
\begin{quote}
\textbf{Definition 13:} (\textit{ordinary least squares}) Let $z_n$ be the natural numbers up to $n$. Further, consider candidate primes $q$ and the candidate factorization $X = F_n(q)$. Then

1) $\ln z_n = F_n(q) \ln q + \varepsilon_n$ is called a "linear factorization model" of $z_n$;

2) $\ln \hat{q} = (X^T X)^{-1} X^T (\ln z_n)$ is called the "log-prime estimator" (LPE);

3) $\ln \hat{q}$ are also referred to as the "estimated log primes";

4) $\ln \hat{z}_n = F_n(q) \ln \hat{q}$ are called the "predicted log integers".
\end{quote}

\begin{quote}
\textbf{Theorem 2:} The LPE minimizes the residual sum of squares.
\end{quote}

\begin{quote}
\textbf{Proof 2:} The residual sum of squares for (\ref{eq:GLM}) is given by

\begin{equation} \label{eq:GLM-RSS1}
\mathrm{RSS}(\beta) = \sum_{i=1}^n \varepsilon_i^2 = \varepsilon^T \varepsilon = (y-X\beta)^T (y-X\beta)
\end{equation}

which can be expanded to

\begin{equation} \label{eq:GLM-RSS2}
\mathrm{RSS}(\beta) = y^T y - y^T X \beta - \beta^T X^T y + \beta^T X^T X \beta
\end{equation}

and differentiated to

\begin{equation} \label{eq:GLM-dRSS}
\mathrm{RSS}'(\beta) = 2 X^T X \beta - 2 X^T y \; .
\end{equation}

Setting this derivative to zero yields

\begin{equation} \label{eq:GLM-OLS}
\hat{\beta} = (X^T X)^{-1} X^T y
\end{equation}

which conforms to the estimated log primes in Def.~13.2. \hspace\fill $\blacksquare$
\end{quote}

With the OLS estimator at hand, we can now consider different cases of candidate primes:

\vspace{0.5em}
\begin{itemize}

\setlength{\itemsep}{0em}
\setlength{\parskip}{0.5em}

\item
Case I: The candidate primes are a real subset of the prime numbers: $q \subset p_n$.

\item
Case II: The candidate primes equal the prime numbers: $q = p_n$.

\item
Case III: The candidate primes are a real superset of the prime numbers: $q \supset p_n$.

\item
Case IV: The candidate primes equal the natural numbers: $q = z_n$.

\end{itemize}
\vspace{0.5em}

Note that these cases generalize the examples from equation (\ref{eq:q-ex}). Figure~2 shows one example for each case and compares (i) the log natural numbers $\ln z_n$ to the predicted log integers $\ln \hat{z}_n$ as well as (ii) the log candidate primes $\ln q$ to the estimated log primes $\ln \hat{q}$, as given in Def.~13. The candidate primes used in the figure are:

\vspace{0.5em}
\begin{itemize}
	
\setlength{\itemsep}{0em}
\setlength{\parskip}{0.5em}

\item
Case I: $q_1 = [3, 5, 11, 17]^T$.

\item
Case II: $q_2 = [2, 3, 5, 7, 11, 13, 17, 19]^T = p_{20}$.

\item
Case III: $q_3 = [2, 3, 4, 5, 7, 8, 11, 12, 13, 15, 17, 19]^T$.

\item
Case IV: $q_4 = [1, 2, 3, \ldots, 18, 19, 20]^T = z_{20}$.
	
\end{itemize}
\vspace{0.5em}

\pagebreak
\begin{center}
\includegraphics[width=1.0\linewidth, clip=true, trim=75 80 45 20]{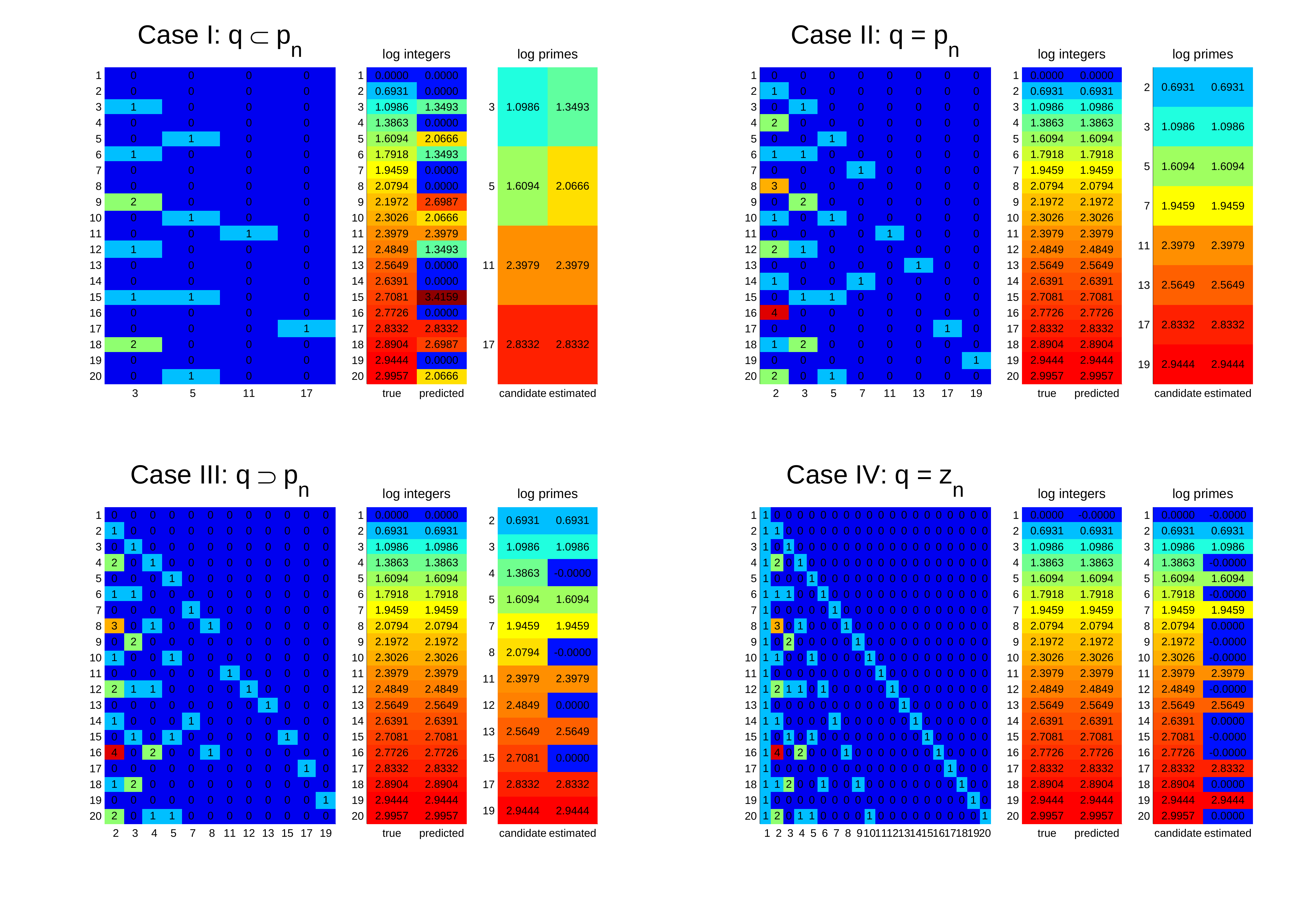}
\end{center}

\textbf{Figure 2.} Four different cases of candidate primes for $n = 20$. Candidate prime vectors are given in the text. Case II is also used in Figure~1 and Case IV is also used in Figure~3. Each panel consists of the candidate factorization (left), comparison of log natural numbers $\ln z_n$ vs. predicted log integers $\ln \hat{z}_n$ (middle) and comparison of log candidate primes $\ln q$ vs. estimated log primes $\ln \hat{q}$ (right). All in all, we make the following observations: (i) As soon as all primes smaller than or equal to $n$ are included in $q$, the log integers are predicted perfectly with maximal accuracy (upper right and lower panels) in which case we call $q$ "complete". If some prime numbers are missing, not all variation in the natural numbers can be captured (upper left panel). (ii) Only if $q = p_n$, the candidate primes are identical to the estimated primes (upper right panel) in which case we call $q$ "valid". If some primes are missing or non-primes are present, there is disagreement (upper left and lower panels). (iii) If $q$ contains all primes smaller than or equal to $n$, non-primes are automatically "switched off" by the LPE and receive a weight of zero whereas primes receive their logarithm as weight (lower panels), consistent with the fundamental theorem of arithmetic (\ref{eq:FTA}) and its logarithmized version (\ref{eq:FTA-ln}). If some elements of $p_n$ are missing in $q$, estimation tends to be unreliable (upper left panel), consistent with the view of the primes as the atomic particles of the natural numbers. (iv) In summary, one can say that the prime numbers are the sparsest set using which one can fully decompose the natural numbers (upper right panel). Equivalently, one could say that the primes are those numbers from the set of all possible candidate primes that minimize the \textit{prediction error} $(\ln z_n-\ln \hat{z}_n)^T (\ln z_n-\ln \hat{z}_n)$ and the \textit{estimation error} $(\ln q-\ln \hat{q})^T (\ln q-\ln \hat{q})$. This refines prime number identification as a model comparison problem in which the least complex from the most accurate models is selected as the optimal solution.

\pagebreak
Based on these observations, we set up consistency conditions for candidate primes and formulate a theorem about the behavior of the LPE for different candidate primes.

\begin{quote}
\textbf{Definition 14:} (\textit{consistency conditions}) Let $q$ be candidate primes and $\ln \hat{q}$ the LPE. Then, we call $q$

1) "valid", if $\ln q = \ln \hat{q}$;

2) "complete", if $\ln z_n = \ln \hat{z}_n$;

3) "consistent", if it is valid and complete.
\end{quote}

\begin{quote}
\textbf{Theorem 3:} Let $q$ be candidate primes and $X = F_n(q)$ the corresponding candidate factorization. Then, $\ln \hat{q} = (X^T X)^{-1} X^T (\ln z_n)$ and:

1) If $q$ is consistent, then $q = p_n$ and vice versa.

2) If $q = p_n$, then $\ln \hat{q} =  \ln p_n$.

3) If $q \supset p_n$, then

\begin{equation} \label{eq:LANT-GLM-superset}
(\ln \hat{q})_j =  \left\{ \begin{matrix} \ln q_j & , & \text{if} \; q_j \in \mathbb{P} \\ 0 & , & \text{if} \; q_j \notin \mathbb{P} \end{matrix} \right., \quad j = 1, \ldots, m \; .
\end{equation}

4) If $q = z_n$, then

\begin{equation} \label{eq:LANT-GLM-naturals}
(\ln \hat{q})_i = \left\{ \begin{matrix} \ln i & , & \text{if} \; i \in \mathbb{P} \\ 0 & , & \text{if} \; i \notin \mathbb{P} \end{matrix} \right., \quad i = 1, \ldots, n \; .
\end{equation}
\end{quote}

\begin{quote}
\textbf{Proof 3:} We prove this theorem step by step.

1) If $q$ is consistent, it follows from Def.~14.1, 14.2 and 13.4 that $\ln z_n = F_n(q) \ln q$. According to Def.~4, there is only one solution for $q$ and this is $q = p_n$. Conversely, if $q = p_n$, then $F_n(q) = F_n$ by Def.~6. We also know that $\ln z_n = F_n \ln p_n$ from Def.~4 which implies that $p_n = q$ is consistent according to Def.~14.3. \hspace\fill $\square$

2) If $q = p_n$, then $F_n(q) = F_n$ by Def.~6. For this case, Def.~4 gives a solution for which $\varepsilon_n = 0_n$, namely $\ln q = \ln p_n$. If there is a solution for which $\mathrm{RSS}(\ln q) = 0$, the LPE must select this solution by Th.~2. Therefore, $\ln \hat{q} = \ln p_n$. \hspace\fill $\square$

3) If $q \supset p_n$, $q$ contains primes and non-primes. For this case, we can construct a solution for which $\varepsilon_n = 0_n$, namely the solution given by (\ref{eq:LANT-GLM-superset}). By Def.~5 and Th.~1, columns of $F_n(q)$ are linearly independent. Therefore, this is the only solution for which $\mathrm{RSS}(\ln q) = 0$. The rest follows the proof of 2). \hspace\fill $\square$

4) This is a special case of 3). \hspace\fill $\square$

This completes the proof. \hspace\fill $\blacksquare$
\end{quote}

\pagebreak
\section{Inverting the Factorization Matrix}

In this section, we want to develop something like a \textit{collective primality test} for the set of all natural numbers up to $n$ with the help of the following theorem:

\begin{quote}
\textbf{Theorem 4:} If $q = z_n$, then $F_n(q)$ is a quadratic $n \times n$ matrix and

1) $F_n(q)$ is invertible;

2) $\ln \hat{q} = \left[ F_n(q) \right]^{-1} (\ln z_n)$.
\end{quote}

\begin{quote}
\textbf{Proof 4:} We prove this theorem step by step.

1) From Def.~6 and Th.~1, it follows that $F_n(z_n)$ is a lower triangular matrix. The determinant of a triangular matrix equals the product of its diagonal entries. Since all diagonal elements of $F_n(z_n)$ are 1, $\det\left[ F_n(z_n) \right] = 1 \neq 0$ and $F_n(z_n)$ is invertible. \hspace\fill $\square$

2) The theory of linear equation systems states that an LES that can be represented as $Ax = b$ with the $n \times n$ invertible matrix $A$ has exactly one solution given by $\hat{x} = A^{-1}b$. Translated to our example, this implies that $\ln \hat{q} = \left[ F_n(z_n) \right]^{-1} (\ln z_n)$. \hspace\fill $\square$

This completes the proof. \hspace\fill $\blacksquare$
\end{quote}

As we now know (i) that the LPE, quite comfortably, "switches off" non-prime entries in a candidate prime vector (see Theorem~3.4) and (ii) that the LPE reduces to a simpler form when the candidate primes equal the natural numbers (see Theorem~4.2), the problem really reduces to finding the inverse of $F_n(z_n)$. We put forward the following solution:

\begin{quote}
\textbf{Theorem 5:} If $q = z_n$, then $F_n(q)$ is a quadratic $n \times n$ matrix and

\begin{equation} \label{eq:Fn-zn-inv}
\left[ F_n(q) \right]^{-1} = \prod_{i=1}^n \left( M_{n+1-i}(2) - \left( \sum_{j=1}^{\lfloor \log_{(n+1-i)}n \rfloor} e_{\overline{(n+1-i)^j}|n} \right) e_{n+1-i|n}^T \right) \; .
\end{equation}
\end{quote}

\begin{quote}
\textbf{Proof 5:} We will use Gauss-Jordan elimination to invert the factorization matrix. This means, we will transform $F_n(z_n)$ into the identity matrix $I_n$ by left-multiplication with elementary matrices and get $[F_n(z_n)]^{-1}$ as the product of these matrices.

Let $F = F_n(z_n)$. According to Th.~1 and Def.~10, the $i$-th column of $F$ is

\begin{equation} \label{eq:Th5-St1}
f_{i|n} = \sum_{j=1}^{\lfloor \log_{i}n \rfloor} e_{\overline{i^j}|n} \; .
\end{equation}

In order to remove $f_{i|n}$ from $F$, i.e. replace it by $e_{i|n}$ to reach $I_n$, we have to left-multiply $F$ with a matrix containing $-f_{i|n}$. However, (i) $-f_{i|n}$ has to be placed into the $i$-th column and (ii) it may not remove the diagonal 1 from $F$. This is achieved by (i) right-multiplying $-f_{i|n}$ with the transposed $e_{i|n}$ and (ii) subtracting it from the $I_n$ that has a 2 in the $i$-th column. In this way, we obtain the extended elementary matrix $E_i$:

\begin{equation} \label{eq:Th5-St2}
E_i = M_i(2) - f_{i|n} \, e_{i|n}^T \; .
\end{equation}

If $F$ is successively left-multiplied with $E_i$, $i = 1,\ldots,n$, it will become $I_n$. Since we have to reverse the order in the product to account for successive left-multiplication, we obtain:

\begin{equation} \label{eq:Th5-St3}
F^{-1} = \prod_{i=1}^n E_{n+1-i} = \prod_{i=1}^n \left( M_{n+1-i}(2) - f_{n+1-i|n} \: e_{n+1-i|n}^T \right) \; .
\end{equation}

With application of (\ref{eq:Th5-St1}), we have

\begin{equation} \label{eq:Th5-St4}
F^{-1} = \prod_{i=1}^n \left( M_{n+1-i}(2) - \left( \sum_{j=1}^{\lfloor \log_{(n+1-i)}n \rfloor} e_{\overline{(n+1-i)^j}|n} \right) e_{n+1-i|n}^T \right)
\end{equation}

which conforms to equation (\ref{eq:Fn-zn-inv}). \hspace\fill $\blacksquare$
\end{quote}

Figure 3 displays an example for inversion of a quadratic candidate factorization matrix.

\begin{center}
\includegraphics[width=1.0\linewidth, clip=true, trim=175 160 125 140]{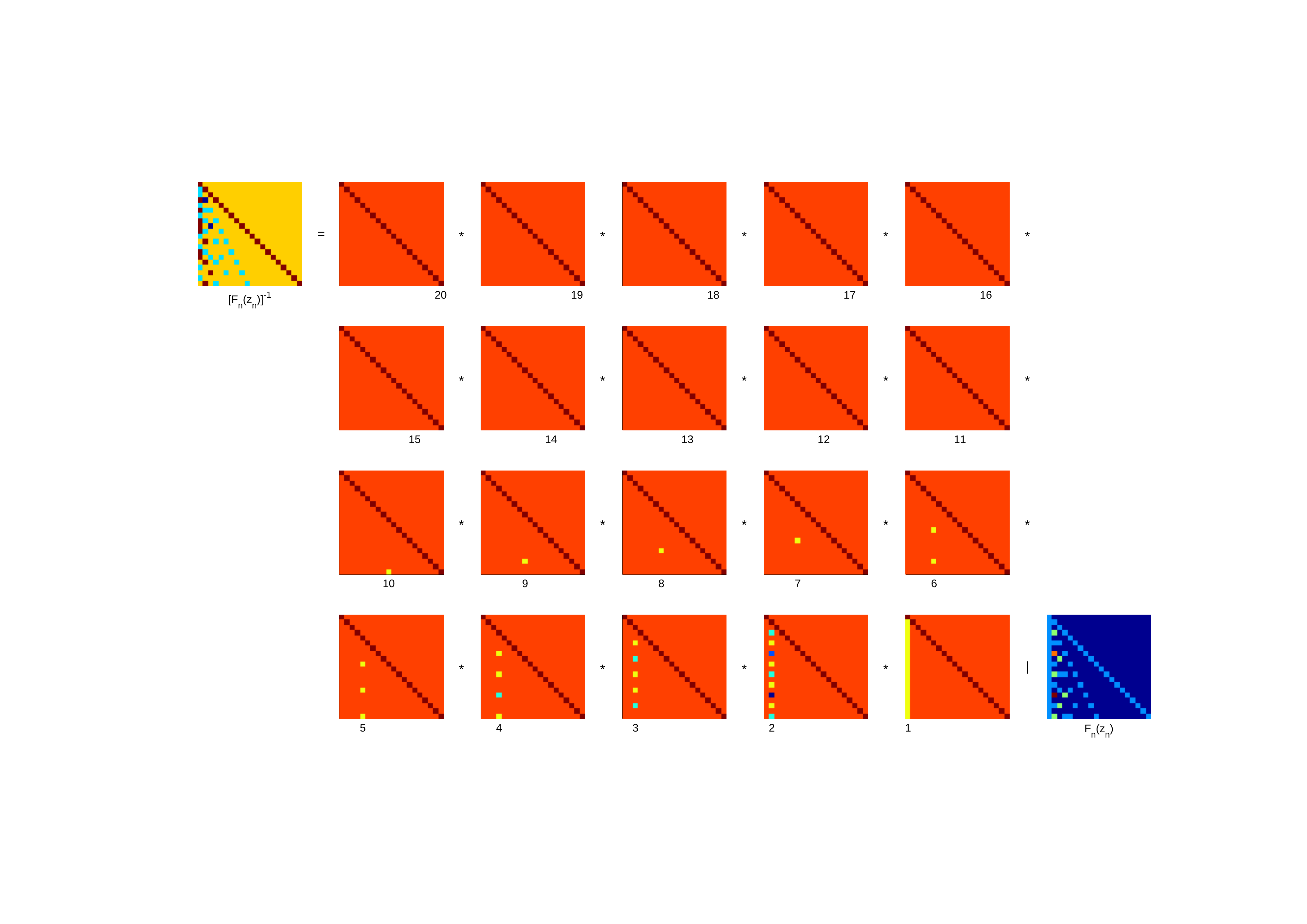}
\end{center}

\textbf{Figure 3.} Example for inverting the factorization matrix with $n = 20$. The candidate factorization $F_n(z_n)$ is shown in the lower right, the extended elementary matrices for left-multiplication $E_n,\ldots,E_1$ are shown in the middle and the resulting matrix inverse $[F_n(z_n)]^{-1}$ is shown in the upper left. Note that the matrices have a different color scale.

\pagebreak
\section{Expressing Number-Theoretic Functions}

In this section, we want to employ the results derived so far, especially Theorem~3, to express certain functions that are important in the field of number theory.

\subsection{IsPrime}

The \textit{IsPrime function} $\mathrm{ip}(x)$ is defined as (OEIS, A010051)

\begin{equation} \label{eq:ip}
\mathrm{ip}(x) = \left\{ \begin{matrix} 1 & , & \text{if} \; x \in \mathbb{P} \\ 0 & , & \text{if} \; x \notin \mathbb{P} \end{matrix} \right. \; .
\end{equation}

Using LANT terminology, $\mathrm{ip}(x)$ can be expressed amazingly simple as

\begin{equation} \label{eq:ip-LANT}
\mathrm{ip}(i) = \frac{(\ln \hat{q})_i}{\ln i} \quad \text{where} \quad \ln \hat{q} = \left[ F_n(z_n) \right]^{-1} (\ln z_n) \quad \text{with} \quad n \geq i
\end{equation}

which is a trivial consequence of Theorem~3.4. Note that $\mathrm{ip}(0)$ is not defined and that $\mathrm{ip}(1)$ is an indeterminate form, consistent with 1 being considered neither prime nor composite. Further, in contrast to $\mathrm{ip}(x)$, $\mathrm{ip}(i)$ is only defined for positive natural numbers.

\subsection{PrimeCount}

The \textit{PrimeCount function} $\pi(x)$ is defined as (OEIS, A000720)

\begin{equation} \label{eq:pi}
\pi(x) = \left| \left\lbrace n \in \mathbb{P} \, | \, n \leq x \right\rbrace \right| \; .
\end{equation}

In the LANT framework, $\pi(x)$ can be expressed similarly simple as

\begin{equation} \label{eq:pi-LANT}
\pi(n) = \sum_{i=2}^n \mathrm{ip}(i) = \sum_{i=2}^n \frac{(\ln \hat{q})_i}{\ln i}
\end{equation}

which follows from equation (\ref{eq:ip-LANT}). Note that this sum starts at $i = 2$, because $\mathrm{ip}(1)$ is an indeterminate form. Again, in contrast to $\pi(x)$ which is defined for any number $x \in \mathbb{R}$ (Platt, 2013), $\pi(n)$ is only defined for positive natural numbers.

\subsection{Chebyshev functions}

The \textit{first Chebyshev function} is given by (Dusart, 2010)

\begin{equation} \label{eq:Ch1}
\vartheta(x) = \sum_{p \leq x} \ln p
\end{equation}

and the \textit{second Chebyshev function} is given by (Dusart, 2010)

\begin{equation} \label{eq:Ch2}
\psi(x) = \sum_{p^k \leq x} \ln p
\end{equation}

where the sums are extending over all prime numbers $p \in \mathbb{P}$ satisfying $p \leq x$ or $p^k \leq x$.

\pagebreak
Again following Theorem~3.4 and based on equation (\ref{eq:ip-LANT}), we have

\begin{equation} \label{eq:Ch1-LANT}
\vartheta(n) = \sum_{i=1}^n (\ln \hat{q})_i = 1_n^T (\ln \hat{q})
\end{equation}

as a simple expression for the first Chebyshev function. The second Chebyshev function cannot be easily represented using LANT quantities, but is related to the first one by

\begin{equation} \label{eq:Ch2-Ch1}
\psi(x) = \sum_{n=1}^\infty \vartheta(x^{1/n}) \; .
\end{equation}

\subsection{von Mangoldt function}

The \textit{von Mangoldt function} is given by (Conrey, 2003)

\begin{equation} \label{eq:vM}
\Lambda(n) = \left\{ \begin{matrix} \ln p & , & \text{if} \; n = p^k, \; p \in \mathbb{P}, \; k \geq 1 \\ 0 & , & \text{otherwise} \end{matrix} \right. \; .
\end{equation}

It can be related to the Chebyshev functions by (Conrey, 2003)

\begin{equation} \label{eq:Ch2-vM}
\psi(x) = \sum_{n \leq x} \Lambda(n) = \sum_{n=1}^{\lfloor x \rfloor} \Lambda(n) \; .
\end{equation}

Using LANT, we obtain the following reformulation of $\Lambda(n)$

\begin{equation} \label{eq:vM-LANT}
\Lambda(i) = \sum_{j=1}^n (\ln \hat{q})_j \left[ \ln i = [F_n(z_n)]_{i,j} \, (\ln \hat{q})_j \right]
\end{equation}

where $[a = b]$ is Iverson bracket notation and $[A]_{i,j}$ refers to the $(i,j)$-th entry of $A$.

\subsection{Riemann $\zeta$ function}

The previously mentioned functions $\pi(x)$, $\vartheta(x)$, $\psi(x)$ and $\Lambda(n)$ are closely related to the complex-valued \textit{Riemann $\zeta$ function} that is given by (Riemann, 1859)

\begin{equation} \label{eq:zeta1}
\zeta(s) = \sum_{n=1}^\infty \frac{1}{n^s} = \frac{1}{1^s} + \frac{1}{2^s} + \frac{1}{3^s} + \ldots 
\end{equation}

which, due to the fundamental theorem of arithmetic, is equivalent to

\begin{equation} \label{eq:zeta2}
\zeta(s) = \prod_{p \in \mathbb{P}} \frac{1}{1-p^{-s}} = \frac{1}{1-2^{-s}} \cdot \frac{1}{1-3^{-s}} \cdot \frac{1}{1-5^{-s}} \cdot \ldots \; .
\end{equation}

Note that these equations only hold for $\mathrm{Re}(s) > 1$, but $\zeta(s)$ can be analytically continued to the the complete real-positive complex half-plane using a Dirichlet eta series by

\begin{equation} \label{eq:zeta3}
\zeta(s) = \frac{1}{1-2^{1-s}} \sum_{n=1}^\infty \frac{(-1)^{n-1}}{n^s} = \frac{1}{1-2^{1-s}} \cdot \left( \frac{1}{1^s} - \frac{1}{2^s} + \frac{1}{3^s} - \ldots \right) \; .
\end{equation}

\pagebreak
\section{The Riemann Hypothesis}

The Riemann $\zeta$ function has trivial zeros at $s = -2, -4, -6, \ldots$ and non-trivial zeros which are known to lie in the critical strip $0 < \mathrm{Re}(s) < 1$. The \textit{Riemann hypothesis} (RH) states that all non-trivial zeros are located on the critical line with real part $\sfrac{1}{2}$:

\begin{equation} \label{eq:RH}
s \in \left\lbrace z \in \mathbb{C} \, | \, \zeta(z) = 0 \wedge \mathrm{Re}(z) > 0 \right\rbrace \Rightarrow \mathrm{Re}(s) = \frac{1}{2} \; .
\end{equation}

RH remains one of number theory's unsolved problems, as it has neither been proven nor falsified so far. However, RH has a lot of important consequences in number theory and is connected to the prime-counting function $\pi(x)$. In particular, it has been shown that RH is equivalent to the following statement (von Koch, 1901):

\begin{equation} \label{eq:RH-eq1}
\pi(x) = \mathrm{Li}(x) + \mathcal{O}(\sqrt{x} \ln x) \; .
\end{equation}

This means that, for a certain $k > 0$ and $x_0 \in \mathbb{R}$, it holds that

\begin{equation} \label{eq:RH-eq2}
|\pi(x) - \mathrm{Li}(x)| \leq k \cdot \sqrt{x} \ln x \quad \text{for all} \quad x \geq x_0 \; .
\end{equation}

Specifically, it has been shown that under RH (Schoenfield, 1976)

\begin{equation} \label{eq:RH-eq3}
|\pi(x) - \mathrm{Li}(x)| \leq \frac{1}{8 \pi} \cdot \sqrt{x} \ln x \quad \text{for all} \quad x \geq 2657 \; .
\end{equation}

In these formulas, $\mathrm{Li}(x)$ is the logarithmic integral function:

\begin{equation} \label{eq:Li}
\mathrm{Li}(x) = \int_{2}^{x} \frac{1}{\ln t} \, \mathrm{d}t \; .
\end{equation}

Remember that we have an explicit formula for $\pi(x)$ (\ref{eq:pi-LANT}) and note how structurally similar this equation is to $\mathrm{Li}(x)$ (\ref{eq:Li}). Therefore, proving RH through means of LANT might be a promising direction. At first sight, there seem to be two strategies:

\vspace{0.5em}
\begin{itemize}

\setlength{\itemsep}{0em}
\setlength{\parskip}{0.5em}

\item
Simplify the left-hand side of (\ref{eq:RH-eq2}) by writing $\pi(x)$ as an integral.

\item
Simplify the left-hand side of (\ref{eq:RH-eq2}) by writing $\mathrm{Li}(x)$ as a sum.

\end{itemize}
\vspace{0.5em}

Incidentally, we have also observed that the function $\mathrm{ld}(n) = \ln \left( \det \left[ F_n^T F_n \right] \right)$ seems to be asymptotically equivalent to $\mathrm{Li}(n)$. This conjecture and other questions will be investigated in future research. A good point to start with might be the further inversion of the factorization matrix (Soch, in prep.).

\pagebreak
\section{References}

\renewcommand{\section}[2]{}


\begin{thebibliography}{9}

\bibitem{Conrey_2003}
Conrey JB (2003): "The Riemann Hypothesis". \textit{Notices of the American Mathematical Society}, vol. 50, no. 3, pp. 341-353.

\bibitem{Dusart_2010}
Dusart P (2010): "Estimates of some functions over primes without R.H.". \textit{arXiv math}, arXiv:1002.0442v1; URL: \url{http://arxiv.org/abs/1002.0442v1}.

\bibitem{Euclid}
Euclid (1996): \textit{Die Elemente. B\"{u}cher I-XIII}. Appeared in \textit{Ostwalds Klassiker der exakten Wissenschaften}, vol. 235. Translated from Greek, edited by Clemens Thaer, with a preface from Wolfgang Trageser, reprinted in 1996, 2\textsuperscript{nd} edition.

\bibitem{Platt_2013}
Platt DJ (2013): "Computing $\pi(x)$ analytically". \textit{arXiv math}, arXiv:1203.5712v3; URL: \url{http://arxiv.org/abs/1203.5712v3}.

\bibitem{Riemann_1859}
Riemann B (1859): "Ueber die Anzahl der Primzahlen unter einer gegebenen Gr\"{o}\ss{}e". \textit{Monatsberichte der Berliner Akademie}, November 1859.

\bibitem{Schoenfield_1976}
Schoenfeld L (1976): "Sharper bounds for the Chebyshev Functions $\theta(x)$ and $\psi(x)$. II". \textit{Mathematics of Computation}, vol. 30, no. 134, pp. 337-360; DOI: 10.2307/2005976.

\bibitem{Soch_in_prep}
Soch J (in prep.): "Linear Algebraic Number Theory, Part II: The Inverse Factorization", in preparation.

\bibitem{von_Koch_1901}
von Koch H (1901): "Sur la distribution des nombres premiers". \textit{Acta Mathematica}, vol. 24, iss. 1, pp. 159-182; DOI: 10.1007/BF02403071.

\bibitem{OEIS}
\textit{The On-Line Encyclopedia of Integer Sequences} (OEIS); URL: \url{http://oeis.org/}.

\end{thebibliography}
\end{document}